\documentclass[twoside,leqno,11pt]{article}

\usepackage{ifthen,amsmath,graphicx,showframe}

\newtheorem{theorem}{Theorem} 

\newtheorem{corollary}{Corollary}

\numberwithin{equation}{section} \numberwithin{lemma}{section}

\numberwithin{theorem}{section} \numberwithin{corollary}{section}

\numberwithin{proposition}{section}

\numberwithin{definition}{section} \numberwithin{example}{section}

\numberwithin{remark}{section}

\def\TITLE       {GEOMETRY AND ITERATION OF DIANALYTIC TRANSFORMATION}
\def\SHORTTITLE  {GEOMETRY AND ITERATION OF DIANALYTIC TRANSFORMATION}
\def\LASTNAMEI   {Cao-Huu}% 

\def\FIRSTNAMEI  {Tuan}%

\def\TOWNI       {Toronto}%

\def\COUNTRYI    {Canada}% 

\def\ADDRESSI    {2275 Bayview Avenue, Toronto, Canada M4N 3M6
}% <--- First Author's Address

\def\EMAILI      {tuan@theory.lcs.mit.edu, tuan@nmr.mgh.harvard.edu} 
\def\LASTNAMEII  {Ghisa}%

\def\FIRSTNAMEII {Dorin}%

\def\TOWNII      {Toronto}% 

\def\COUNTRYII   {Canada}% <--- 2nd Author's Country

\def\ADDRESSII   {2275 Bayview Avenue, Toronto, Canada M4N 3M6
}% <--- 2nd Author's Address

\def\EMAILII     {dghisa@yorku.ca}% <--- 2nd Author's Email Address

\def\LASTNAMEIII {}% <--- 3rd Author's Last Name

\def\FIRSTNAMEIII{}% <--- 3rd Author's First Name

\def\TOWNIII     {}% <--- 3rd Author's Residence Town

\def\COUNTRYIII  {}% <--- 3rd Author's Country

\def\ADDRESSIII  {}% <--- 3rd Author's Address

\def\EMAILIII    {}% <--- 3rd Author's Email Address

\def\CLASSIFICATION{AUTOMATION, COMPUTERS, APPLIED MATHEMATICS\\
Volume 13, Number 1, 2004,
ISBN 1221-437X\\
MSC 2000: 30F50, 37F50}% <--- Optional

\def\KEYWORDS{Nonorientable Klein surface, dianalytic transformation, 
Julia set, Fatou set, dianalytic dynamics, Blaschke product
}% <-- Optional

 {} % <---

\def\ABSTRACT{{\bf ABSTRACT --\/}
This paper is a continuation of the paper  [5] dealing
 with dynamics
of dianalytic transformations of nonorientable Klein surfaces. We are
examining mainly the transformations of the real projective plane $P^{2}, $
whose orientable double cover is the Riemann sphere $\overline{C}$ . It is
shown that the automorphisms of \ $P^{2}$ are projections of the rotations of
\ $\overline{C}$ and some of the other dianalytic transformations of \ $P^{2}$
are projections of Blaschke products.

}% <-- Optional

\begin{document}

%\frontmatter++++++++++++++++++++++++++++++++++++++++++++++++++++%

\begin{center}

\par\TITLE

\par(\SHORTTITLE)

\par\FIRSTNAMEI \ \LASTNAMEI

\par\TOWNI \  \COUNTRYI

\par\ADDRESSI

\par\EMAILI

\par\medskip

\par\FIRSTNAMEII \ \LASTNAMEII

\par\TOWNII \ \COUNTRYII

\par\ADDRESSII

\par\EMAILII

\par\medskip

\par\FIRSTNAMEIII \ \LASTNAMEIII

\par\TOWNIII \ \COUNTRYIII

\par\ADDRESSIII

\par\EMAILIII

\par\CLASSIFICATION{}

\par\KEYWORDS{}

\par\ABSTRACT{}

\end{center}

%++++++++++++++++++++++++++++++++++++++++++++++++++++++++++++++++++++++++++

%\frontmater
\newenvironment{proof}[1][Proof]{\textbf{#1.} }{\ \rule{0.5em}{0.5em}}

\section{Introduction}

  It is known that the only analytic maps $F:\overline{C}->\overline{C} $
 are the rational functions.
  If $F(z)=p(z)/q(z),$ where \ $p$ \ and \ $q$ \ are polynomials with
complex coefficients that have no common factors, then the degree of $F$ \ is
by definition,
     $\deg(F)=\max\{\deg(p),\deg(q)\}$,
and it represents the number (counted with multiplicity) of inverse images of
any point \ $z\in\overline{C}$ \ by the map $F.$

  The Fatou-Julia theory applies to rational maps $F$ \ whose degree is
at least two. The dynamics of M\oe bius transformations is in general simpler
and in most cases can be completely described.

  The use of computer graphics, as well as the application of
quasiconformal mapping techniques has brought in the last decades an
extraordinary development \ of Fatou-Julia theory. Moreover, parallel studies
have been done on non rational analytic functions and extensions have been
considered to Riemann surfaces, multivalued functions, several variables, etc.
Our aim is to deal with dynamics on nonorientable Klein surfaces.

  Such a study was initiated in [5], where a complete list of dianalytic
self maps of the real projective plane, the pointed real projective plane and
the Klein bottle has been given. The dynamics of these transformations can be
described by using their lifting to the respective orientable double covers.
In this paper we take a closer look to these dynamics.

  Let us consider on the Riemann sphere the conformal structure\ (see
[2]) induced by the atlas
     $\Upsilon=\{(C,\varphi_{1}),(\overline
{C}-\{0\},\varphi_{2})\},$
where \ $\varphi_{1}(z)=z,$ \ $\varphi_{2}(z)=1/z.$
  The antianalytic involution \ $h:\overline{C}->\overline{C}$ \ defined
by \ $h(z)=-1/\overline{z}$ \ if \ $z\neq0,$ \ $z\neq\infty,$ \ and
$h(0)=\infty,$ \ $h(\infty)=0$ \ makes \ $(\overline{C},h)$ \ into a symmetric
Riemann surface. The real projective plane \ $P^{2}=\overline{C\text{ }}/<h>$
\ endowed with the unique dianalytic structure which makes the canonical
projection \ $\pi:\overline{C}->P^{2}$ \ a dianalytic function, represents a
nonorientable Klein surface (see [3] and [4]).

\section{Dianalytic Transformations of \ $P^{2}$}

We have shown in [5] that any dianalytic transformation \ $f$ \ of \ $P^{2}$
\ has the representation
  $f(\widetilde{z})=\widetilde{F(z)}$, \ \ or
\ \ $f(\widetilde{z})=\widetilde{F(\overline{z}})$,
where \ $\widetilde{z}=(z,h(z)\}$ \ and \ $F$ \ is a rational function of the form: 

(1) \ \ \ \ \ \ \ \ \  \ \ \ $F(z)=e^{i\theta}\frac{a_{0}z^{2n+1}+a_{1}%
z^{2n}+...+a_{2n+1}}{-\overline{a}_{2n+1}z^{2n+1}+\overline{a}_{2n}%
z^{2n}-...+\overline{a}_{0}},$ \ \ \ \ \ $|a_{0}|+|a_{2n+1}|\neq0,$
$\ \theta\in R.$

Also, any dianalytic automorphism of \ $P^{2}$ \ has the representation
   $g(\widetilde{z})=\widetilde{G(z)},$ \ \ or
\ \ $g(\widetilde{z})=\widetilde{G(\overline{z}}),$ \ \ \ where

(2)\ \ \  \ \ \ \ \ \ \ \ \ $G(z)=e^{i\theta}\frac{az+b}{-\overline{b}z+\overline{a}},$
\ \ $|a|+|b|\neq0,$ \ \ $\theta\in R$

\noindent By using these facts we can prove the following theorem.

\begin{theorem} Any dianalytic transformation \ $f$ \ of \ $P^{2}$ \ has the representation
     $f(\widetilde{z})=\widetilde{F(z)}$ \ \ \ or
\ \ \ $f(\widetilde{z})=\widetilde{F(\overline{z}}),$ \ \ where

(3)   \ \ \  \ \ \ \ \ \ \ \ \ $F(z)=e^{i\alpha}\prod_{k=1}^{2n+1}\frac{z-z_{k}%
}{1+\overline{z_{k}}z},$ \ \ $\alpha\in R,$ \\ for some arbitrary , not
necessarily distinctive

\ $\ \ \ \ \ \ \ \ \ \ \ \ \ \ \ \ \ \ \ \ \ \ \ \ \ \ \ \ z_{k}\in C.$
\end{theorem}
\begin{proof} \ Indeed, at least one of the polynomials in \ (1) \ has the degree
\ $2n+1$. \ Suppose that this is the polynomial at the numerator. Then
\ $a_{0} $ $\neq0,$ \ and the polynomial has the decomposition
   $a_{0}(z-z_{1})(z-z_{2})...(z-z_{2n+1})$,
where the numbers \ $z_{k}$ \ are repeated, if they are multiple roots.

\noindent  By
making the change of variable \ $z=-1/\overline{w}$ \ at the denominator, this becomes:

 $   \frac{\overline{a_{2n+1}}}{\overline{w}^{2n+1}%
}+\frac{\overline{a_{2n}}}{\overline{w}^{2n}}+...+\overline{a_{0}}%
=\frac{1}{\overline{w}^{2n+1}}(\overline{a_{0}}\overline{w}^{2n+1}%
+\overline{a_{1}}\overline{w}^{2n}+...+\overline{a_{2n+1}})=$

    $\frac{1}{\overline{w}^{2n+1}}(\overline{a_{0}%
w^{2n+1}+a_{1}w^{2n}+...+a_{2n+1}})=\frac{1}{\overline{w}^{2n+1}}%
\overline{a_{0}\prod_{k=1}^{2n+1}(w-z_{k})}=$

    $\overline{a_{0}}\prod_{k=1}^{2n+1}\frac{\overline
{w-z_{k}}}{\overline{w}^{2n+1}}=\overline{a_{0}}\prod_{k=1}^{2n+1}%
\frac{\overline{w}-\overline{z_{k}}}{\overline{w}^{2n+1}}=\overline{a_{0}%
}\prod_{k=1}^{2n+1}(1-\frac{\overline{z_{k}}}{\overline{w}})=\overline{a_{0}%
}\prod_{k=1}^{2n+1}(1+\overline{z_{k}}z)$

\noindent
Then \ $F$ \ has the expression:

    $F(z)=e^{i\theta}\frac{a_{0}\prod_{k=1}^{2n+1}%
(z-z_{k})}{\overline{a_{0}}\prod_{k=1}^{2n+1}(1+\overline{z_{k}}%
z)}=e^{i(\theta+2\beta)}\prod_{k=1}^{2n+1}\frac{z-z_{k}}{1+\overline{z_{k}}%
z}=e^{i\alpha}\prod_{k=1}^{2n+1}\frac{z-z_{k}}{1+\overline{z_{k}}z}$

where \ \ $\beta=\arg a_{0}$ \ and \ \ $\alpha=\theta+2\beta.$
\end{proof}

  An analogous reasoning brings us to the same result if we suppose that
the polynomial at the denominator has the degree \ $2n+1.$

\section{Dynamics of Dianalytic Automorphisms of \ $P^{2}$}

Theorem 1 has the following Corollary:
\begin{corollary}
 The dianalytic automorphisms of \ $P^{2}$ \ have the representation: 

     g($\widetilde{z})=\widetilde{G(z)}$ \ \ \ \ or
\ \ \ \ $g(\widetilde{z})=\widetilde{G(\overline{z}}),$ 
where \ $G$ \ is a rotation of the Riemann sphere.
\end{corollary}
\begin{proof} We need to show that \ $z->$ \ $e^{i\alpha}\frac{z-z_{0}}%
{1+\overline{z_{0}}z}$ \ \ is indeed a rotation of the Riemann sphere. It is
obvious that \ $w->e^{i\alpha}w$ \ \ represents a rotation of the Riemann
sphere around the vertical diameter. On the other hand, it can be easily
checked that \ $z->\frac{z-z_{0}}{1+\overline{z_{0}}z}$ \ has the fixed points
\ $ie^{i\theta}$ \ and \ $-ie^{i\theta},$ \ where \ $\theta=\arg z_{0},$
\ therefore this mapping represents a rotation of the Riemann sphere around
the diameter passing through the two fixed points. The composition of the two
rotations give a new rotation of the Riemann sphere around a third diameter.
The extremities of that diameter have as stereographic projections the fixed
points of the transformation \ $G.$ The dynamics of \ $G$ \ is then described
by the Steiner's net associated to those fixed points (see, for example [1],
page 85). It can be easily shown that \ if \ $z_{0}=re^{i\theta},$ while the
rotation axis of \ \ $z->\frac{z-z_{0}}{1+\overline{z_{0}}z}$\ depends on
\ $\theta,$ the rotation angle is\ $\ \cos^{-1}\frac{r^{2}-1}{r^{2}+1}.$
Therefore, whether the final rotation is rational or irrational depends on
\ $\alpha$ \ as well as on \ $r.$ \ In the case of rational rotations
we have cyclic iterations. Due to the $h-$invariance of \ $G$, the order of
the cycle should be even. Then \ $g$ \ will have also cyclic iterations of an
order half of that of \ $G$.
\end{proof}

\section{Dianalytic Maps by Blaschke Products}

A Blaschke product is a function of the form 

\ (4)\ \ \ \ \ \ \ \ \ \ \ \ \ \ \ \ \ \ \ \ \ \ \ \ \ \ \ \ \ \ $F(z)=e^{i\theta
}\prod\frac{z-z_{k}}{1-\overline{z}_{k}z},$ \ \ $|z_{k}|<1$

with a finite or infinite number of terms. When there are infinitely many
terms, the product is written usually under the form:

\ \ \ \ \ \ \ \ \ \ \ \ \ \ \ \ \ \ \ \ \ \ \ \ \ \ \ \ \ \ \ \ \ \ \ \ \ \ \ \ $B(z)=$%
\ \ $\prod\frac{\overline{z}_{k}}{|z_{k}|}\frac{z_{k}-z}{1-\overline{z}_{k}z}$

and it is known that if

(5)$\ \ \ \ \ \ \ \ \ \ \ \ \ \ \ \ \ \ \ \ \ \ \ \ \ \ \ \ \ \ \ \ \ \ \ \ \sum
(1-|z_{k}|)<\infty,$

then it converges uniformly on compact subsets of \ $\overline{C}-E,$ where
\ $E$ $\subseteq\partial D$\ is the set of accumulation points of \ $\{z_{k} $
$\}$ (here $\partial D$ is the unit circle).

We notice that \ $F$ \ is not generally \ $h-$invariant, therefore it is not
of the form (3). However, if for every \ $k$ \ there is \ $k^{\prime}$ \ such
that \ $z_{k}=-z_{k^{\prime}}$ , then: 

\ \ \ \ \ \ \ \ \ \ \ \ \ \ \ \ \ \ \ \ \ \ \ \ \ $\frac{z-z_{k}}%
{1-\overline{z}_{k}z}\cdot\frac{z-z_{k^{\prime}}}{1-\overline{z}_{k^{\prime}%
}z}=\frac{z-z_{k}}{1-\overline{z}_{k^{\prime}}z}\cdot\frac{z-z_{k^{\prime}}%
}{1-\overline{z}_{k}z}=\frac{z-z_{k}}{1+\overline{z}_{k}z}\cdot
\frac{z-z_{k^{\prime}}}{1+\overline{z}_{k^{\prime}}z}$

is \ $h-$invariant. If such a condition is fulfilled and for an odd number of
\ $k$, \ $z_{k}=0,$\ and if the product is finite, then \ (4) \ is of the form
\ (3) \ and the functions \ $f$ \ defined by

\ \ \ \ \ \ \ \ \ \ \ \ \ \ \ \ \ \ \ \ \ \ \ \ \ \ \ \ \ \ \ \ \ \ \ $\widetilde
{z}->\widetilde{F(z)},$ \ respectively \ \ $\widetilde{z}->\widetilde
{F(\overline{z})}$

are dianalytic transformations of \ $P^{2}.$ Then 

(6)\ \ \ \ \ \ \ \ \ \ \ \ \ \ \ \ \ \ \ \ \ \ \ \ \ \ \ \ \ \ \ \ \ \ \ \ \ \ \ \ \ $F(z)=e^{i\theta
}z^{2p+1}\prod_{k=1}^{m}\frac{z^{2}-z_{k}^{2}}{1-\overline{z}_{k}^{2}z^{2}},$
\ $z_{k}\neq0,$ \ $k=1,2,...,m$

It is obvious that for such a function the unit disk remains invariant. We
ignore the trivial case where \ $m=0$ \ and suppose that \ $\deg(F)\geq3.$ By
the Schwarz Lemma , for any \ $z,$ with \ $|z|<1,$ \ we have \ $|F(z)|<|z|.$
Consequently, the sequence \ $(|F^{\circ n}(z)|)$ \ is decreasing, which shows
that the open unit disk $D$ \ belongs to the Fatou set of \ $F$. \ Moreover,
$D$ \ is the basin of attraction of the attracting fixed point \ $0.$ By a
known result, (see [6], or [8], p. 46), the Julia set of $\ F$ \ must be the
whole unit circle. Consequently, we have proved:

\begin{theorem} If \ $f:P^{2}->P^{2}$ is the projection of a finite
Blaschke product \ $F$ \ with \ $\deg(F)\geq3$ then the Fatou set of \ $f$
\ is \ $\{\widetilde{z}:|z|<1\}$ \ and the Julia set of \ $f$ \ is
\ $\{\widetilde{z}:|z|=1\}.$ 
\end{theorem}

It follows that the dynamics of those dianalytic self-maps of \ $P^{2}$, which
happen to be projections of finite Blaschke products, are very simple. In [8],
p. 150, more general finite Blaschke products were considered, namely for
which some of \ $z_{k}$ \ are inside the unit circle and some are outside. It
has been shown that such a product might possess Herman rings, therefore the
structure of the Julia set and that of the Fatou set is much more complicated.
For an \ $h-$invariant Blaschke product of that type, the respective Herman
rings should be symmetric, such that they project in pairs on the same Herman
ring in \ $P^{2}.$

Suppose now that the infinite Blaschke product \ $B(z)$ \ is \ $h-$invariant
and the convergence condition (5) is fulfilled. This last condition implies
that \ $B(z)$ \ is analytic in \ $\overline{C}$\ \ with the exception the
poles \ $1/\overline{z}_{k}$ \ and the points of \ $E,$ \ which are essential
singularities (see [1]). We have again that \ $|B(z)|<1$ \ for every \ $z$
$\ $with $\ |z|<1$ $\ $and \ $|B(z)|>1$ \ for every \ $z$ \ with \ $|z|>1.$
Then, in this case too, the unit disk is completely invariant (see [6]), and
so is the exterior of the unit disk.

%\bigskip\noindent
\section{References}

\noindent
[1]. \ Ahlfors, L.V. , Complex Analysis, McGraw-Hill Book Company, 1979

\noindent
[2]. \ Ahlfors, L.V. \& \ Sario, L., Riemann Surfaces, Princeton University
Press, 1960

\noindent
[3]. \ Alling, N \& Greenleaf, N., The Foundation of Theory of Klein Surfaces,
Lecture Notes in
     Mathematics, No. 219, \ Springer Verlag, 1971

\noindent
[4]. \ Andreian Cazacu, C., On the Morphisms of Klein Surfaces, Rev. Roumaine
Math. Pures
    \ Appl. 31, 1986, 461-470\

\noindent
[5]. Barza, I. \& Ghisa, D., Dynamics of Dianalytic Transformations of Klein 
Surfaces,
    \ Mathematica Bohemica, No. 2, 2004, 129-140\

\noindent
[6]. \ Blanchard, P., Complex Analytic Dynamics on the Riemann Sphere, Bull.
of Am. Math.
    \ Soc., Vol. 11, No. 1, 1984, 85-141\ \

\noindent
[7]. Lehto, O., Univalent functions and Teichm\"{u}ller Spaces,
Springer-Verlag, 1987

\noindent
[8]. Milnor, J., Dynamics in One Complex Variable, Friedr. Vieweg \& Sohn,
1999

\end{document}